\documentclass[a4paper,oneside,11pt]{scrartcl}
\input{mypaper.sty}
\usepackage{enumerate}
\usepackage{subfigure}
\usepackage{tikz}
\usetikzlibrary{calc,angles,quotes,decorations.markings}
\usepackage{cleveref}

\title{On low-dimensional approximation of function spaces of interior regularity}
\author{
	S. Aziz\footnote{Faculty of Mathematics, Physics and Computer Science, University of Bayreuth, 95447 Bayreuth, Germany}
	\and
	M. Bauer$^*$
	\and
	M. Bebendorf$^*$
	\and
	T. Rau\footnote{UL Solutions SIS, Erlangen, Germany}
}
\date{\today}

\begin{document}
\maketitle

\begin{abstract}
Many elliptic boundary value problems exhibit an interior regularity property, which can be exploited to construct local approximation spaces that converge exponentially within function spaces satisfying this property. 
These spaces can be used to define local ansatz spaces within the framework of generalised finite element methods, leading to a better relation between dimensionality and convergence order.
In this paper, we present a new technique for the construction of such spaces for Lipschitz domains. Instead of
the commonly used approach based on eigenvalue problems it relies on extensions of approximations performed on the boundary.
Hence, it improves the influence of the spatial dimension on the exponential convergence and allows to construct the local spaces by solving the original kind of variational problems on easily structured domains. 
\end{abstract}
%\textbf{Keywords:} ACA,  BEM, non-local operators, hierarchical matrices, fast solvers

\section{Introduction}
We consider the efficient numerical solution of variational problems
\begin{equation}\label{eq:bvp}
\text{find } u\in V:\quad a(u,\fie)=\ell(\fie),\quad \fie\in C_0^\infty(\Omega),
\end{equation}
with a bilinear form $a$ and a Lipschitz domain $\Omega\subset\R^d$. 
If $a$ results from the variational formulation of a second order elliptic boundary value problem, the solution $u$ of such problems, and its discrete approximation from finite element spaces, can usually be approximated using methods of linear or logarithmic-linear complexity. Some of the most prominent examples are \textit{multigrid methods}~\cite{BrambleMG,MG,MR1193013}, \textit{hierarchical matrices}~\cite{MR2000c:65039,MR2001i:65053}, and \textit{$hp$-finite element methods}~\cite{doi:10.1137/0718033,JMelenk02}. While multigrid methods can be regarded as iterative methods that exploit smoothness with respect to a sequence of nested grids, $hp$-methods rely on a combination of grid refinement and the local approximation from polynomial spaces of relatively high degree. The aim of this article is to show that the general concept of
\textit{interior regularity} can be used to devise new numerical methods of logarithmic-linear complexity. Since interior regularity does not require $u$ to be globally smooth and since in addition to linear elliptic problems this property is observed also for nonlinear problems, we expect this new approach to be equally applicable to a broad range of problems
such as problems with non-smooth coefficients and boundary value problems for the $p$-Laplacian.    

The approach presented in this article shares the idea of approximation of local solution spaces
\[
X_{\ell,g}(D):=\{u\in V: a(u,\fie)=\ell(\fie)\text{ for all }\fie\in C_0^\infty(D)\text{ and }u=g\text{ on }\partial D\cap\partial \Omega\},
\]
where $D\subset\Omega$, with \textit{generalized finite element methods} (GFEM)~\cite{IBCaOs94,MMIB96}.
The latter is constructed by partitioning the computational domain~$\Omega$ into a set of subsets $D$ and constructing finite-dimensional approximation spaces over each subset for the local solution space~$X_{\ell,g}(D)$.
GFEM combines local spaces via a partition of unity, which allows to treat the local problems
in parallel prior to the solution of the global problem, which typically has a significantly smaller number of degrees of freedom than usual finite element methods. The accuracy of the GFEM solution is controlled by the local approximation error \cite{MMIB96}.
 In the case of the  \textit{multiscale spectral generalized finite element methods} (MS-GFEM)~\cite{IBRL11,CMRSTD22} the local approximation
spaces are constructed in an optimal way using the solution of suitable eigenvalue problems.
Since the numerical treatment of eigenvalue problems is usually quite costly, the aim of this article is to
construct local approximation spaces via a recursive approximation technique that shares some principles with multigrid methods, i.e., the method presented in this article is based on a sequence of local variational problems of the original
type~\eqref{eq:bvp}. If for $D$ easily structured domains such as disks/balls or rectangles/boxes are used, then
highly optimized solvers such as multigrid can be employed for the solution of the local problems.
The constructed spaces will be seen to converge exponentially to the local solution space~$X_{\ell,g}(D)$.
Exponential convergence can also be observed in to $hp$-FEM if the solution is sufficiently smooth so that it can be
approximated by polynomials. In the approach of this article higher regularity will not be required. 
The presented method will rely on the minimum assumption that is used for the convergence of finite element methods, i.e.\ throughout this article it will be used that $u\in H^{1+\alpha}(\Omega)$ with some arbitrarily small but fixed $\alpha>0$. 
It is also worth mentioning that the technique used to generate the basis of the proposed exponentially convergent approximation spaces resembles the technique used in the \textit{virtual element method} (VEM)~\cite{VEM1BP,VEM2HG}, in the sense that the degrees of freedom are primarily located on the boundary of the domain.

In \cite{MBWHInv} we have constructed approximation spaces which converge exponentially
w.r.t.\ the $L^2$-norm for the approximation of harmonic functions
\[
X_{0,0}(D):=\{u\in H^1(D): a(u,\fie)=0\text{ for all }\fie\in C_0^\infty(D)\text{ and }u=0\text{ on }\partial D\cap\partial \Omega\}
\]
in the case of convex domains~$D$.
The convergence proof is not constructive as it relies on the existence of an $L^2$-projection onto the (closed) space~$X_{0,0}(D)$.
While the $L^2$-norm is suitable in the context of hierarchical matrix approximations (cf.~\cite{Bebendorf:2008,WH15}), the $H^1$-norm is the natural choice in the context of solutions of second order boundary value problems. Hence, the aim of this article is to generalise and improve these results in several directions:
\begin{enumerate}[(i)]
	\item estimates with respect to arbitrary Sobolev norms will be presented,
	\item general Lipschitz domains will be considered,
	\item a constructive approach will be presented,
	\item the dependence of the dimension of the approximation space $\Xi_\eps$ on the spatial dimension~$d$ will be improved.
\end{enumerate}
Although the construction presented in Sect.~\ref{sec:one} relies on general interior regularity estimates,
in Sect.~\ref{sec:two} we confine ourselves to Dirichlet boundary value problems and bilinear forms $a(u,v)=\int_\Omega \nabla v^TC \nabla u\ud x$ with a symmetric positive definite matrix $C(x)\in\R^{d\times d}$ for almost all $x\in\Omega$. 
 While for $L^2$-estimates the $H^1$-regularity of solutions is sufficient,
for $H^1$-estimates a regularity higher than $H^1$ is required. Since Lipschitz domains will be considered throughout this article, the regularity of solutions is typically $H^{1+\alpha}$ with some $0<\alpha\leq1$ if the coefficients $c_{ij}$ are sufficiently smooth; see~\cite{Ne64}. This requires proving interior regularity estimates for Sobolev norms of fractional order. 
The approach presented in~\cite{Borm:2010aa} to $H^1$-estimates avoids fractional order estimates but the convergence proof still uses the (non-constructive) projection onto~$X_{0,0}(D)$. In order to overcome this, a method to implement the projection onto $X_{0,0}(D)$ is required.
Our new approach uses the harmonic extension of approximations constructed on the boundary, which also requires
$H^{1+\alpha}$-regularity.

Sect.~\ref{sec:three} compares several numerical techniques for the implementation of the harmonic extension in the case of
Dirichlet problems for the Laplace equation. The extension
is either done using Green's functions or the method of fundamental solutions for both the disk and rectangular domains.
In the case of disks we also compare our construction with trigonometric approximation.
Furthermore, problems with the construction of $\Xi_\eps$ using finite element approximations of the
harmonic extension are reported.

While in this article we focus on the construction of local approximation spaces, the next step is to use these spaces
in combination with a suitable covering of the computational domain $\Omega$ by easily structured subsets such as disks/balls or rectangles/boxes for the construction of a generalized finite element method.

\section{General Setting}\label{sec:one}
Let $\Omega\subset\R^d$ denote the computational domain and $D\subset\Omega$ be a bounded Lipschitz domain.
For $s\in\R_0^+$ let
\[
H^s(D)=\{v\in H^k(D):\partial^\beta v\in H^\sigma(D)\text{ for all }\beta\in\N_0^d,\,|\beta|= k\}
\]
denote the Sobolev-Slobodeckij space of fractional order $s=k+\sigma$, $k\in\N_0$, $\sigma\in(0,1)$, wherein we define
$H^\sigma(D)=\{v\in L^2(D):|v|_{H^\sigma(D)}<\infty\}$ 
with the Sobolev-Slobodeckij semi-norm \[|v|_{H^\sigma(D)}^2=\int_D\int_D \frac{|v(x)-v(y)|^2}{|x-y|^{d+2\sigma}}\ud y\ud x.\]
$H^s(D)$ is a Hilbert space equipped 
with the norm $\norm{v}_{H^s(D)}=\left(\norm{v}_{H^k(D)}^2+|v|_{H^s(D)}^2\right)^{1/2}$ 
and the semi-norm \[|v|_{H^s(D)}:=\left(\sum_{|\beta|=k}|\partial^\beta v|_{H^\sigma(D)}^2\right)^{1/2}.\]
Additionally, we define the space $H^s_0(D)$ as the completion of $C_0^\infty(D)$ w.r.t.\ $\norm{\cdot}_{H^s(D)}$.
For $s\neq1/2$ the dual space of $H_0^s(D)$ is denoted with $H^{-s}(D)$. 

We consider linear spaces of functions $X(D)\subset H^s(D)$ such that for their restriction to an open set $K\subset D$ having positive distance $\dist(K,\Omega\cap\partial D)$ to the boundary of $D$ within $\Omega$
it holds that $X(D)|_K\subset X(K)$. The elements $u\in X(D)$ are assumed to satisfy an interior regularity estimate 
\begin{equation}\label{eq:innreg}
\norm{u}_{H^r(K)}\leq \frac{c_R}{\dist^{r-s}(K,\Omega\cap\partial D)}\norm{u}_{H^s(D)}
\end{equation}
with a real number $r>s$ and a constant $c_R>0$.
We assume that the higher regularity in $K$ can be exploited in the sense that
for all $n\in\N$ a linear space $V_n(K)\subset X(K)$
with $\dim{V_n(K)}\leq n$ exists such that for every $u\in X(D)$ there
is $v_u\in V_n(K)$ such that
\begin{equation}\label{eq:exlwralg}
  \norm{u-v_u}_{H^s(K)}\leq c_A\left(\frac{\diam\,K}{\sqrt[m]{n}}\right)^{r-s}\norm{u}_{H^r(K)}
\end{equation}
with some $c_A\geq 1$ and $m\geq1$.
If the distance of $K$ to the boundary of $D$ within $\Omega$ is large relative to its diameter, i.e.
\begin{equation}\label{eq:adm}
  \eta\,\dist(K,\Omega\cap\partial D)\geq \diam\, K
  \end{equation}
with some parameter $\eta>0$, then the algebraic convergence w.r.t.\ the dimension of the
approximating space assumed in~\eqref{eq:exlwralg} can be improved to an exponential one using a recursive construction similar to the technique presented in \cite{MBWHInv}.

\begin{theorem}
Assuming~\eqref{eq:innreg}--\eqref{eq:adm}, for every $\eps\in (0,1)$ there is a subspace $\Xi_\eps(K)\subset X(K)$ with $\dim \Xi_\eps(K)\lesssim |\log\eps|^{m+1}$  such that for all $u\in X(D)$ there is $\xi_u\in \Xi_\eps(K)$
satisfying  \label{thm:one}
\[
\norm{u-\xi_u}_{H^s(K)}\leq \eps\norm{u}_{H^s(D)}.
\]
\end{theorem}
\begin{proof}
We consider a sequence of $L:=\lceil |\log\eps|/(r-s)\rceil$ sets $K_L=K$ and
$K_j=\{x\in\Omega:\dist(x,K)<\rho_j\}$ with $\rho_j:=(1-j/L)\,\dist(K,\Omega\cap\partial D)$ for $j=0,\dots,L-1$.
Notice that $K_L\subset K_{L-1}\subset\dots\subset K_0\subset D$.

According to \eqref{eq:exlwralg} there are subspaces $V_n(K_j)\subset X(K_j)$ with
\[\dim V_n(K_j)\leq n:=\lceil (c_Ac_R\eps^{-1/L})^{1/(r-s)}(2+\eta)\rceil^m L^m\]
and for all $w\in X(K_{j-1})$ there is an element $v_w\in V_n(K_j)$ such that 
\[
\norm{w-v_w}_{H^s(K_j)}\leq c_A\left(\frac{\diam\,K_j}{\sqrt[m]{n}}\right)^{r-s}\norm{w}_{H^r(K_{j})}.
\]
We apply this approximation recursively to $r_0:=u\in X(K_0)$. Setting 
\begin{equation}\label{eq:recapp}
r_j:=r_{j-1}|_{K_j}-v_j,\quad j=1,2,\dots,L,
\end{equation}
where $v_j\in V_n(K_j)$ denotes the approximation of the restriction of~$r_{j-1}\in X(K_{j-1})$ to $K_j$, we obtain that $r_j\in X(K_j)$ as $X(K_{j-1})|_{K_j}\subset X(K_j)$
and
\[
\norm{r_j}_{H^s(K_j)}\leq c_A\left(\frac{\diam\,K_j}{\sqrt[m]{n}}\right)^{r-s}\norm{r_{j-1}}_{H^r(K_j)}.
\]
Using \eqref{eq:innreg}, this leads to
\[
\norm{r_j}_{H^s(K_j)}\leq c_Ac_R\left(\frac{\diam\,K_j}{\sqrt[m]{n}\,\dist(K_j,\Omega\cap\partial K_{j-1})}\right)^{r-s}\norm{r_{j-1}}_{H^s(K_{j-1})}.
\]
Since $\dist(K_j,\Omega\cap\partial K_{j-1})=\rho_{j-1}-\rho_j=\dist(K,\Omega\cap\partial D)/L$, it follows from \eqref{eq:adm} that
\[\diam\, K_j\leq \diam\, K+2\rho_j\leq (2+\eta)\,\dist(K,\Omega\cap\partial D)=(2+\eta)L\,\dist(K_j,\Omega\cap\partial K_{j-1}).\]
Hence, using the definition of $n$
\[
\norm{r_j}_{H^s(K_j)}\leq c_Ac_R\left(\frac{(2+\eta) L}{\sqrt[m]{n}}\right)^{r-s}\norm{r_{j-1}}_{H^s(K_{j-1})}
\leq \eps^{1/L}\,\norm{r_{j-1}}_{H^s(K_{j-1})}.
\]
The recursive application of the previous estimate for $j=L,\dots,1$ yields
\[\norm{r_L}_{H^s(K_L)}\leq \eps^{1/L} \norm{r_{L-1}}_{H^s(K_{L-1})}\leq\dots\leq \eps\norm{r_0}_{H^s(K_0)}\leq \eps\norm{u}_{H^s(D)}.\]
Notice that $r_L=u|_K-\xi_u$ with \[\xi_u:=\sum_{j=1}^L v_j|_K\in \Xi_\eps(K):=\bigoplus_{j=1}^L V_n(K_j)|_K\subset X(K)\] and $\dim \Xi_\eps(K)\leq L n\leq \lceil (c_A c_R)^{1/(r-s)}\e(2+\eta)\rceil^m L^{m+1}$
due to $\eps^{-1/L}\leq \e^{r-s}$.
\end{proof}

In some sense, \eqref{eq:recapp} shares principles with multigrid procedures. While in multigrid methods smoothing is required to make the error amenable to approximation on a coarser grid, here the restriction to subdomains increases the smoothness of the remainder.

\section{Application to Laplace-Type Problems}\label{sec:two}
In \cite{MBWHInv} we have used the previous technique for constructing finite-dimensional approximation spaces which converge exponentially w.r.t.\ the $L^2$-norm for the approximation of harmonic functions
\[
X_{0,0}(D):=\{u\in H^1(D): a(u,\fie)=0\text{ for all }\fie\in C_0^\infty(D)\text{ and }u=0\text{ on }\partial D\cap\partial \Omega\}.
\]
There and in what follows we consider the bilinear form
\begin{equation}\label{eq:varf}
	a(u,v)=\int_D \nabla v^T C \nabla u\ud x
\end{equation}
with $C(x)\in\R^{d\times d}$ being symmetric positive definite for almost all $x\in D$, i.e.\ there is $\lambda,\Lambda>0$ such
that
\[
\lambda\norm{v}_2^2\leq v^TC(x)v\leq \Lambda\norm{v}_2^2,\quad v\in\R^d,\;x\in D.
\]

While the $L^2$-norm is suitable in the context of hierarchical matrix approximations (cf.~\cite{Bebendorf:2008,WH15}), the $H^1$-norm is the natural choice in the context of solutions of boundary value problems, i.e., w.r.t.\ the $H^1$-norm we consider the approximation  of affine spaces
\[
X_{\ell,g}(D):=\{u\in H^1(D): a(u,\fie)=\ell(\fie)\text{ for all }\fie\in C_0^\infty(D)\text{ and }u=g\text{ on }\partial D\cap\partial \Omega\}.
\]
$X_{\ell,g}(D)$ is the space of local
solutions of the variational form
\[
\text{find }u\in \{v\in H^1(\Omega):v=g\text{ on }\partial\Omega\} \text{ such that }a(u,v)=\ell(v)\text{ for all }v\in C_0^\infty(\Omega)
\]
 of the Dirichlet problem for the Poisson equation
with general linear forms $\ell\in H^{-1}(\Omega)$ and $g\in H^{1/2}(\partial\Omega)$.
While $D$ was assumed to be convex in \cite{MBWHInv}, in this article we consider general Lipschitz domains~$D$.

The lack of smoothness of the domain requires to use fractional Sobolev spaces. For Lipschitz domains it is known (see~\cite{Ne64}) that elements of~$X_{0,0}(D)$ have a regularity slightly higher than $H^1(D)$, i.e.\ $X_{0,0}(D)\subset H^{1+\alpha}(D)$ with some $0<\alpha\leq1$ if the coefficients $c_{ij}$ are sufficiently smooth.
We rely on the boundedness
\begin{equation}\label{eq:conttr}
	\norm{\gamma v}_{H^{s-1/2}(\partial D)}\lesssim \norm{v}_{H^s(D)},\quad v\in H^s(D),
\end{equation}
of the trace operator $\gamma:H^s(D)\to H^{s-1/2}(\partial D)$ for $1/2<s<3/2$; see~\cite[Thm.~3.38]{McLean}.

In order to show an interior regularity estimate of type~\eqref{eq:innreg},
we first prove a Poincar\'e-Friedrichs inequality for Sobolev spaces of fractional order using a well-known compactness argument.

\begin{lemma}
 For $\sigma\in(1/2,1)$ there is a constant $c_{\sigma,D}>0$ such that 
  \[
\norm{u}_{L^2(D)}\leq c_{\sigma,D}|u|_{H^\sigma(D)}
  \]
  holds for all $u\in H^\sigma_0(D)$.
\end{lemma}
\begin{proof}
  Suppose that there is no constant $c>0$
  such that $\norm{u}_{H^\sigma(D)}\leq c\,|u|_{H^\sigma(D)}$ for all $u\in H_0^\sigma(D)$.
Then there is a sequence $\{u_n\}_{n\in\N}\subset H_0^\sigma(D)$ with $\|u_n\|_{H^\sigma(D)}>n\,|u_n|_{H^\sigma(D)}$,
and $v_n:=u_n/\|u_n\|_{H^\sigma(D)}$ is a bounded sequence in $H_0^\sigma(D)$. Due to the compact embedding $H^\sigma(D)\hookrightarrow L^2(D)$, 
there is a subsequence $\{v_{n_k}\}_{k\in\N}$ which converges in~$L^2(D)$.
In particular, $\{v_{n_k}\}_{k\in\N}$ is a Cauchy sequence in~$L^2(D)$.
Since 
\[
|v_n|_{H^\sigma(D)}=\frac{|u_n|_{H^\sigma(D)}}{\norm{u_n}_{H^\sigma(D)}}<\frac{1}{n},  
\]
we see that $\lim_{k\to\infty} |v_{n_k}|_{H^\sigma(D)}=0$. From
\[
\norm{v_{n_i}-v_{n_j}}_{H^\sigma(D)}^2\leq \norm{v_{n_i}-v_{n_j}}^2_{L^2(D)}+\left[|v_{n_i}|_{H^\sigma(D)}+|v_{n_j}|_{H^\sigma(D)}\right]^2,
\]
it follows that $\{v_{n_k}\}_{k\in\N}$ is a Cauchy sequence in the complete space $H_0^\sigma(D)$, which converges to $v\in H_0^\sigma(D)$.
From $|v|_{H^\sigma(D)}=\lim_{k\to\infty}|v_{n_k}|_{H^\sigma(D)}=0$ and the definition of $|\cdot|_{H^\sigma(D)}$ we obtain that $v\in H_0^\sigma(D)$ is constant and thus $v=0$. This shows the contradiction as on the other hand
$\norm{v}_{H^\sigma(D)}=\lim_{k\to\infty}\norm{v_{n_k}}_{H^\sigma(D)}=1$.
\end{proof}

The constant $c_{\sigma,D}>0$ in the previous lemma depends on~$D$ but is not known in general.
Using a scaling argument, we see that it has to be proportional to the $\sigma$-th power of the diameter of~$D$, i.e.
\begin{equation}\label{eq:PFs}
  \norm{u}_{L^2(D)}\leq \hat c_{\sigma,D} (\diam\,D)^\sigma \,|u|_{H^\sigma(D)}.
\end{equation}

Next, we prove an interior regularity estimate of type~\eqref{eq:innreg}. It is known (see~\cite{Ev97}) that
$u\in H^2(K)$ provided that the diffusion coefficients satisfy $c_{ij}\in C^1(D)$. However, we will neither require nor be able to benefit from estimates with respect to the $H^2$-norm as the
trace operator in the case of Lipschitz domains is continuous only for functions in $H^{1+\alpha}(K)$ with $1+\alpha<3/2$.

\begin{theorem}
 Let $K \subset D$ be a Lipschitz domain satisfying \eqref{eq:adm} and let $\alpha<1/2$. If the coefficient matrix~$C$ in \eqref{eq:varf} is assumed to consist of entries~$c_{ij}\in C^1(D)$, then there is $c_\alpha>0$ such that
  \[
  \norm{u}_{H^{1+\alpha}(K)} \leq  \frac{c_\alpha}{\dist^\alpha(K,\Omega\cap\partial D)} \norm{u}_{H^1(D)}
  \]
  for all $u \in X_{0,0}(D)$. \label{thm2}
  \end{theorem}
  
  \begin{proof}
  For $K \subset D$  there exists a cut-off function $\chi \in C^\infty(D)$ with $0\leq \chi \leq 1$, $\chi=0$ on $\partial D\cap\Omega$, $\chi= 1$ in $K$, and $\norm{\partial^\beta \chi}_{L^\infty(D)} \leq c/\dist^{|\beta|}(K,\Omega\cap\partial D)$ for all $\beta\in \N_0^d$; see \cite[Thm.~3.6]{McLean}. Using the Fourier transform it can be seen that
   the bilinear form $a$ is well-defined on $H^{1+\alpha}(D)\times H^{1-\alpha}(D)$; see~\cite[Lem.~3.1]{Necas}. Hence, for $u\in X_{0,0}(D)\subset H^{1+\alpha}(D)$ and $v \in H_0^{1-\alpha}(D)$ we define the linear functional 
  \[
  \delta(v) := a(\chi u_0, v) - a(u_0, \chi v),\quad u_0:=u-\overline{u},
  \]
  where $\overline{u}:=\frac{1}{|D|}\int_D u\ud x$ denotes the average of $u$ on $D$.
  Integration by parts yields
  \begin{align*}
  \delta(v) 
  &= \int_D \nabla v^TC\nabla(\chi u_0) -\nabla(\chi v)^TC\nabla u_0\ud x\\
  &= \int_D u_0 \nabla v^TC \nabla\chi + \chi \nabla v^TC \nabla u_0 -v\nabla \chi^TC \nabla u_0 - \chi\nabla v^TC \nabla u_0\ud x  \\
  &= \int_D u_0 \nabla v^TC \nabla \chi  - v \nabla \chi^TC \nabla u_0\ud x 
  = - \int_D v\,\div (u_0 C\nabla \chi)   +v\nabla \chi^TC\nabla u_0 \ud x \\
  &= -2 \int_D v\nabla \chi^TC \nabla u_0\ud x  - \int_D u_0v\,\div(C\nabla\chi)\ud x
  \end{align*}
  and thus
  \begin{align*}
  |\delta(v)| &\leq 2 \norm{C\nabla u_0}_{L^2(D)} \norm{\nabla \chi}_{L^\infty(D)} \norm{v}_{L^2(D)} +
  \norm{\div(C\nabla \chi)}_{L^\infty(D)} \norm{u_0}_{L^2(D)} \norm{v}_{L^2(D)} \\
  &\leq  c\left( 2   + \frac{\diam\,D}{\dist(K,\Omega\cap\partial D)}  \right)  \frac{\norm{\nabla u}_{L^2(D)}}{\dist(K, \Omega\cap\partial D)} \norm{v}_{L^2(D)}.
  \end{align*}
The last estimate follows from the Poincar\'e inequality and $\nabla u_0=\nabla u$.
  Applying \eqref{eq:PFs} to $v\in H_0^{1-\alpha}(D)$
  shows due to \eqref{eq:adm} and $\diam\,D\leq \diam\,K+ 2\,\dist(K,\Omega\cap\partial D)\leq (2+\eta)\,\dist(K,\Omega\cap\partial D)$ that
  \begin{align*}
  |\delta(v)| 
  &\leq c\,\hat c_{1-\alpha,D}\, (4+\eta) \frac{(\diam\,D)^{1-\alpha}}{\dist(K,\Omega\cap\partial D)} \norm{\nabla u}_{L^2(D)} \norm{v}_{H^{1-\alpha}(D)}\\
  &\leq c\,\hat c_{1-\alpha,D}\,(4+\eta) \frac{(2+\eta)^{1-\alpha}}{\dist^\alpha(K,\Omega\cap\partial D)} \norm{\nabla u}_{L^2(D)} \norm{v}_{H^{1-\alpha}(D)}
  \end{align*}
  and thus 
  \begin{equation}\label{eq:dK}
    \norm{\delta}_{H^{\alpha-1}(D)}=\sup_{0\neq v\in H^{1-\alpha}_0(D)} \frac{|\delta(v)|}{\norm{v}_{H^{1-\alpha}(D)}}\leq  \frac{\tilde c_\alpha}{\dist^\alpha(K,\Omega\cap\partial D)} \norm{\nabla u}_{L^2(D)}
  \end{equation}
  with $\tilde c_\alpha:=c\,\hat c_{1-\alpha,D}\,(4+\eta)(2+\eta)^{1-\alpha}$.
Hence, $\delta\in H^{\alpha-1}(D)\subset H^{-1}(D)$. Let $\tilde{u}\in H_0^1(D)$ be the unique solution of the variational problem
  \begin{equation*}
  a(\tilde{u},v) = \delta(v)\quad\text{for all }v\in H^1_0(D).
  \end{equation*}
 The regularity result for Lipschitz domains \cite{Ne64} yields
  \begin{equation*}
  |\tilde{u}|_{H^{1+\alpha}(D)} \lesssim \norm{\delta}_{H^{\alpha-1}(D)}. 
  \end{equation*} 
  Since $u_0\in X_{0,0}(D)$ we have $a(\chi u_0, v) = a(u_0, \chi v) + \delta(v) = \delta(v)$ for all $v\in H^1_0(D)$.
Since $\chi u_0\in H_0^1(D)$, the uniqueness of~$\tilde u$ implies $\tilde{u} = \chi u_0$ and therefore
  \[
  |u|_{H^{1+\alpha}(K)}=|u_0|_{H^{1+\alpha}(K)} \leq |\chi u_0|_{H^{1+\alpha}(D)} \lesssim \norm{\delta}_{H^{\alpha-1}(D)}. 
  \]
  Applying \eqref{eq:dK} and using $\norm{u}_{H^{1+\alpha}(D)}^2=\norm{u}_{H^1(D)}^2+|u|_{H^{1+\alpha}(D)}^2$ together with $\dist(K,\Omega\cap\partial D)\leq\diam\,D$ leads to the assertion.
  \end{proof}
  
  \begin{remark}
  The constants depend on the diffusion coefficient $C$. With a more sophisticated technique it can be proved (cf.~\cite{MB14}) that the contrast, i.e.\ the ratio of the largest and smallest eigenvalue, of $C$ enters the dimension of the approximation spaces only logarithmically.
  \end{remark}

As obviously $X_{0,0}(D)|_K\subset X_{0,0}(K)$, it remains to show the existence of a finite-dimensional approximation space $V_n(K)\subset X_{0,0}(K)$ satisfying~\eqref{eq:exlwralg}.
Let $\mathcal{T}_h$ be a quasi-uniform boundary mesh of~$\partial K$ with $n$ vertices.
For the approximation space $\mathcal{S}_h^1(\partial K)$ of piecewise linear, globally continuous functions
on~$\mathcal{T}_h$ the approximation property
\begin{equation}\label{eq:apprh}
\min_{v_h\in \mathcal{S}^1_h(\partial K)} \norm{u-v_h}_{H^s(\partial K)}\leq ch^{r-s}|u|_{H^r(\partial K)}  
\end{equation}
holds for $u\in H^r(\partial K)$ and $s\in[0,1]$, $r\in[s,1/2+\alpha]$; see \cite[Thm.~10.9]{Steinb07}. 
The space~$V_n(K)$ is defined as the harmonic extension of $\mathcal{S}_h^1(\partial K)$ to $K$
\begin{equation}\label{eq:Vn}
V_n(K):=\{v\in X_{0,0}(K)\text{ such that }\gamma v\in \mathcal{S}_h^1(\partial K)\}. 
\end{equation}

\begin{lemma}
Let $K$ be a Lipschitz domain and $\alpha<1/2$. The previous construction yields a linear space $V_n(K)\subset X_{0,0}(K)$ of dimension $\dim V_n(K)=n$ such that for all $u\in X_{0,0}(K)$ there is $v_u\in V_n(K)$ satisfying \label{lem2}
\[
\norm{u-v_u}_{H^1(K)}\leq c_A\left(\frac{\diam\,K}{\sqrt[d-1]{n}}\right)^\alpha \norm{u}_{H^{1+\alpha}(K)}.
\]
\end{lemma}
\begin{proof}
  Let $u\in X_{0,0}(K)$ be given. As $u\in H^{1+\alpha}(K)$, we apply \eqref{eq:apprh} to its trace $\gamma u\in H^{1/2+\alpha}(\partial K)$.
This yields $v_h\in \mathcal{S}_h^1(\partial K)$ such that
  \[
  \norm{u-v_h}_{H^{1/2}(\partial K)}\leq ch^\alpha|u|_{H^{1/2+\alpha}(\partial K)}.
  \]
With $|\partial K|\sim (\diam\,K)^{d-1}$, the previous estimate can be expressed in terms of the number of degrees of freedom~$n$ as 
  \[
  \norm{u-v_h}_{H^{1/2}(\partial K)}\lesssim \left(\frac{\diam\,K}{\sqrt[d-1]{n}}\right)^\alpha|u|_{H^{1/2+\alpha}(\partial K)}.
  \]
  Let $v_u\in H^1(K)$ be the solution of the
  variational problem $a(v_u,\fie)=0$ for all $\fie\in C_0^\infty(K)$ such that $\gamma v_u=v_h$ on $\partial K$.
  Then $a(u-v_u,\fie)=0$ for all $\fie\in C_0^\infty(K)$ and $\gamma(u-v_u)=\gamma u-v_h$ on $\partial K$.
  Hence, $v_u\in V_n(K)$ and the Lax-Milgram theorem shows
  \[
 \norm{u-v_u}_{H^1(K)}\lesssim \norm{u-v_h}_{H^{1/2}(\partial K)}.
  \]
  Using the boundedness \eqref{eq:conttr} of the trace operator, we obtain
  \[
    \norm{u-v_u}_{H^1(K)} \lesssim \left(\frac{\diam\,K}{\sqrt[d-1]{n}}\right)^\alpha|u|_{H^{1/2+\alpha}(\partial K)}
    \lesssim \left(\frac{\diam\,K}{\sqrt[d-1]{n}}\right)^\alpha\norm{u}_{H^{1+\alpha}(K)}.
  \]
\end{proof}

From the previous proof it can be seen that a slightly higher regularity than $H^1(K)$ is required in order to benefit from
the reduced dimension when discretizing the boundary.

The following corollary proves the exponentially convergent approximation of~$X_{\ell,g}(D)$.
\begin{corollary}
Let $K\subset D$ be a Lipschitz domain satisfying \eqref{eq:adm}. If the coefficient matrix~$C$ in~\eqref{eq:varf} is assumed to consist of entries~$c_{ij}\in C^1(D)$, then for every $\eps>0$ there is an affine subspace $\Xi_\eps(K)\subset X_{\ell,g}(K)$ with $\dim \Xi_\eps(K)\lesssim |\log\eps|^d$ such that for all $u\in X_{\ell,g}(D)$ there is $\xi_u\in \Xi_\eps(K)$ with
\[
\norm{u-\xi_u}_{H^1(K)}\leq \eps \norm{u}_{H^1(D)}. 
\]
\end{corollary}
\begin{proof}
Let $u_{\ell,g}\in X_{\ell,g}(D)$ be defined by $u_{\ell,g}=0$ on $\partial D\cap\Omega$. Then the Lax-Milgram theorem and the boundedness of the trace operator show $\norm{u-u_{\ell,g}}_{H^1(D)}\lesssim
\norm{u}_{H^{1/2}(\partial D)}\lesssim \norm{u}_{H^1(D)}$ for $u\in X_{\ell,g}(D)$. Notice that $u-u_{\ell,g}\in X_{0,0}(D)$.
Setting $m=d-1$, $s=1$ and $r=1+\alpha$ with $\alpha<1/2$, Theorem~\ref{thm2} yields~\eqref{eq:innreg} with $c_R=c_\alpha$. Furthermore, Lemma~\ref{lem2} shows that there is $V_n(K)\subset X_{0,0}(K)$ such that
$\dim V_n(K)=n$ and
\[
\min_{v\in V_n(K)}\norm{w-v}_{H^1(K)}\leq c_A\left(\frac{\diam\,K}{\sqrt[d-1]{n}}\right)^\alpha \norm{w}_{H^{1+\alpha}(K)},  \quad w\in X_{0,0}(K).
\]
Theorem~\ref{thm:one} shows the existence of $\tilde\Xi_\eps(K)\subset X_{0,0}(K)$ with $\dim \tilde\Xi_\eps(K)\leq c \lceil|\log\eps|\rceil^d$, where $c=\lceil \alpha^{d/(1-d)} (c_Ac_\alpha)^{1/\alpha}\e(2+\eta)\rceil^{d-1}$ such that
\[
\min_{\tilde\xi\in \tilde\Xi_\eps(K)}\norm{w-\tilde\xi}_{H^1(K)}\leq \eps \norm{w}_{H^1(D)},\quad w\in X_{0,0}(D).
\]
Defining $\Xi_\eps(K)=u_{\ell,g}|_K+\tilde\Xi_\eps(K)$, we have
\[
\min_{\xi\in\Xi_\eps(K)}\norm{u-\xi}_{H^1(K)}=\min_{\tilde\xi\in\tilde\Xi_\eps(K)}\norm{u-u_{\ell,g}-\tilde\xi}_{H^1(K)}\leq \eps \norm{u-u_{\ell,g}}_{H^1(D)}\lesssim \eps\norm{u}_{H^1(D)}.
\]
\end{proof}

\begin{remark}
	Since the type of the domain $D$ will be fixed (ball, box, etc.), finding $u_{\ell,g}$ in a practical implementation can be done in advance. Notice that the regularity parameter $\alpha$ enters the dimension estimate for $\Xi_\eps(K)$ only via the constant in front of the logarithm.
\end{remark}

The following corollary presents an (improved) $L^2$-estimate. Notice that the construction of $V_n(K)$ via the boundary allows to reduce the dimensionality, i.e.\ the exponent $m$ in the dimension estimate of the approximation space~$\Xi_\eps(K)$ compared to the construction presented in \cite{MBWHInv}.
\begin{corollary} \label{cor:erreps}
	Let the assumptions of the previous corollary be satisfied. Then for every $\eps>0$ there is an affine subspace $\Xi_\eps(K)\subset X_{\ell,g}(K)$ with $\dim \Xi_\eps(K)\lesssim |\log\eps|^d$ such that for all $u\in X_{\ell,g}(D)$ there is $\xi_u\in \Xi_\eps(K)$ with
	\[
	\norm{u-\xi_u}_{L^2(K)}\leq \eps \norm{u}_{L^2(D)}.
	\]
\end{corollary}
\begin{proof}
	We apply the previous corollary to $K$ and $D':=\{x\in D:\dist(x,K)<\frac{1}{2}\dist(K,\Omega\cap\partial D)\}$.
	The constructed set $\Xi_\eps(K)\subset X_{\ell,g}(K)$ satisfies
	\[
	\min_{\xi\in \Xi_\eps(K)} \norm{u-\xi}^2_{L^2(K)}\leq \eps\norm{u}^2_{H^1(D')}\leq \eps\left(\norm{u}_{L^2(D)}^2+\norm{\nabla u}^2_{L^2(D')}\right).
	\]
	The Caccioppoli inequality shows that
	\[
	\norm{\nabla u}_{L^2(D')}\leq \frac{c}{\dist(D',\Omega\cap\partial D)}\norm{u}_{L^2(D)}=\frac{2c}{\dist(K,\Omega\cap\partial D)}\norm{u}_{L^2(D)}.
	\]
	\end{proof}

\section{Numerical Methods} \label{sec:three}
In the following, we present numerical experiments to verify the error estimates from \Cref{cor:erreps} for smooth and Lipschitz domains, respectively. In these we use the bilinear form $a(u,v)=\int_D\nabla u\cdot\nabla v\ud x$.

\subsection{Smooth Domains}
Let $K \subset \R^2$ be a disk with radius $a$. The first step is to generate the basis of the space $\mathcal{S}^1_h(\partial K)$. To do so, polar coordinates $(r,\theta)$ are used and the boundary $\partial K$ is discretised in terms of $\theta$.
Let $\mathcal{T}_h$ be a boundary mesh on $\partial K$ defined as
\begin{equation*}
    \mathcal{T}_h := \left \{ \left [ \theta_j,\theta_{j+1} \right] : 0\leq j \leq n-1 \right \} ,
\end{equation*}
where $0=\theta_0 < \ldots < \theta_n = 2\pi$, with uniform mesh size $h:=2\pi a/n$. It is clear that $\mathcal{T}_h$ is periodic, meaning $\theta_0 = \theta_n$, so the mesh has exactly $n$ intervals and $n$ distinct nodes.

\begin{figure}[ht]
	\centering
	\includegraphics[width=0.48\linewidth]{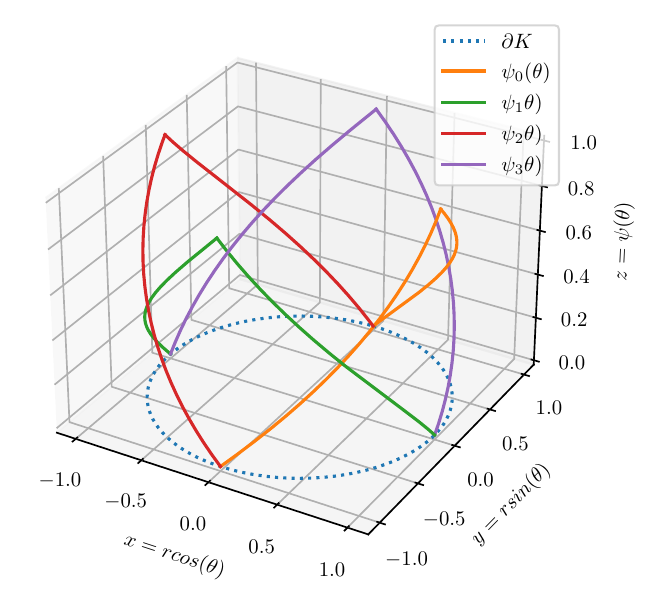}
    \caption{Piecewise linear functions $\psi_i$ on $\partial K$ for $\omega=\pi/2$. \label{fig:sd_hf}}	
\end{figure}

Let $\omega:=h/a$ be the discretisation parameter of~$\theta$. We then define the basis functions $\psi_i$ as
\begin{equation}
\psi_i(\theta) = 
  \begin{cases}
    \frac{\theta-(i-1)\omega}{\omega}, & (i-1)\omega \leq \theta \le i\omega, \\
    \frac{(i+1)\omega -\theta}{\omega}, & i\omega \leq \theta \le (i+1)\omega,\\
    0, &\textnormal{otherwise.}
  \end{cases}
  \label{eq:sd_hat_func}
\end{equation}
The shapes of the basis functions $\psi_i$ for different discretisation parameters~$\omega$ are shown in Figure~\ref{fig:sd_hf}.

The next step is to construct the basis functions $\phi_i$, $i = 1,\ldots,n$, of the space $V_n(K)$ defined in~\eqref{eq:Vn}, which are the harmonic extension of the boundary functions $\psi_i$, i.e., $\phi_i$ are the unique solutions of the following boundary value problems
\begin{subequations}     \label{eq:sd_harm_bas}
    \begin{alignat}{3}
        -\Delta \phi_i(r,\theta) &= 0 &\quad& \text{in } K ,\\
        \phi_i(a,\theta) &= \psi_i(\theta) &\quad& \text{on } \partial K.
    \end{alignat}	
\end{subequations}
The boundary value problems \eqref{eq:sd_harm_bas} can be solved using various numerical schemes. For our experiments, we consider Green's functions method~\cite{sneddon} and the method of fundamental solutions~\cite{Karageorghis}.

\subsection*{Green's functions method}
The Green's functions method is a numerical technique used to solve Dirichlet boundary value problems; see~\cite{sneddon}. It relies on the use of so-called Green's functions $\mathcal{G}$ that satisfy
\begin{alignat*}{3}
    -\Delta_y \op G(x,y) &= \delta_0(y-x), &\quad& y \in K,\\
    \op G(x,y) &= 0, &\quad&  y \in \partial K,
\end{alignat*}
for each $x$ in $K$. Here, $\delta_0$ denotes the Dirac distribution. Consequently the solutions $\phi_i$ of \eqref{eq:sd_harm_bas} at a point $x=(r,\theta)$ in $K$ can be represented as
\[
    \phi_i(r,\theta) = - \int_{\partial K} \partial_{\nu_y} \mathcal{G}(x,y)\, \psi_i(y) \ud s_y,
\]
where $\nu_y$ denotes the unit normal vector at $y\in\partial K$.
For a disk with radius $a$, the Green's function is given by
	\begin{equation*}
		\mathcal{G}(x,y) = - \left [S(x-y) - S\left(|x| \left(\frac{y}{a} -\frac{a x}{|x|^2}\right)\right)\right]
	\end{equation*}
with the singularity function $S(x):=\frac{1}{2\pi}\log|x|$. In this case $\phi_i$ becomes
	\begin{equation*}
		\phi_i(r,\theta) = \frac{1}{2\pi}\int_0^{2\pi} \frac{a^2-r^2}{a^2+r^2-2ar\cos(\theta-\Tilde{\theta})} \psi_i(\Tilde{\theta}) \ud\Tilde{\theta}.
	\end{equation*}
With the definition \eqref{eq:sd_hat_func} of $\psi_i$, the basis functions $\phi_i$ of $V_n(K)$ are then given by
\begin{equation*}
    \phi_i(r,\theta) = \frac{1}{2\pi\omega} \int_{(i-1)\omega}^{i\omega} \frac{(a^2-r^2)(\Tilde{\theta}-(i-1)\omega)}{a^2+r^2-2ar\cos(\theta-\Tilde{\theta})} \ud\Tilde{\theta} + \int_{i\omega}^{(i+1)\omega} \frac{(a^2-r^2)((i+1)\omega - \Tilde{\theta})}{a^2+r^2-2ar\cos(\theta-\Tilde{\theta})} \ud\Tilde{\theta}.
\end{equation*}
The integral above has no general closed-form solution and is only defined within the domain's interior. This is because it exhibits a singularity at $r=a$ and $\theta = \Tilde{\theta}$, which complicates the calculation of the solution near the boundary. Therefore, Gauss quadrature rules tailored to the singularity or adaptive quadrature rules, such as Gauss-Kronrod rules~\cite{kronrod}, are needed to evaluate the integral. 

The main drawback of Green's functions method is that it is difficult to calculate the solution near the boundary and the quality of the solution degrades as the boundary is approached. The following section explores an alternative method for calculating the basis functions.

\subsection*{Method of Fundamental Solutions}
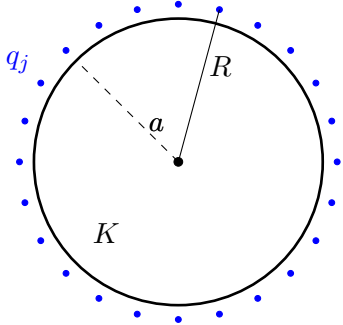
\begin{figure}[ht]
	\centering
\begin{tikzpicture}[scale=1.9]
    \draw (0,0) circle (1);
    \node at (120:0.3) {$a$};
    \draw[line width=1pt] (0,0) circle (1);
    \node at (-0.5,-0.5) {$K$};
    \draw[dashed] (0,0) -- (135:1);
    \node at (120:0.3) {$a$};
    \node at (150:1.1) [above left,blue] {$q_j$};
    \foreach \i in {0,15,...,345}
    \filldraw [blue] (\i:1.1) circle (0.02);
    \draw (0,0) -- (75:1.1) node[midway, above right] {$R$};
    \filldraw  (0,0) circle (0.03);
\end{tikzpicture}
\caption{Singularities placed outside a circular domain.\label{fig:sd_mfs_const}}
\end{figure}

The Method of Fundamental Solutions (MFS) is a meshless collocation boundary method that solves certain elliptic boundary value problems; see~\cite{Karageorghis}. Given a second-order elliptic operator and its fundamental solution~$S$, the MFS represents the solution of the boundary value problem as a linear combination of $N$ fundamental solutions with singularities $q_j$, $j=1,\dots,N$, positioned outside the domain $\overline{K}$
\begin{equation*}
    \tilde \phi_i(x) := \sum_{j=1}^N c_j^{(i)} S(x-q_j)\approx \phi_i(x), \quad x \in \overline{K}.
\end{equation*}
Discretising the boundary of $K$ into $M$ collocation points $x_k\in\partial K$, $k=1,\dots,M$, and applying the Dirichlet boundary conditions
\begin{equation*}
    \tilde\phi_i(x_k) = \sum_{j=1}^N c_j^{(i)} S(x_k-q_j) = \psi_i(x_k),\quad k=1,\dots,M,
\end{equation*}
results in a system of linear equations that can then be solved by the least squares method.\\
For a disk with radius $a$, singularities are positioned on a circle of radius $R>a$, as shown in Figure~\ref{fig:sd_mfs_const}. The number of collocation points $M$ is set at a higher resolution than the singularities, typically $M\approx3N$. The number and positions of the singularities are heuristically determined.

\begin{figure}[ht]
     \centering
        \includegraphics[width=0.48\linewidth]{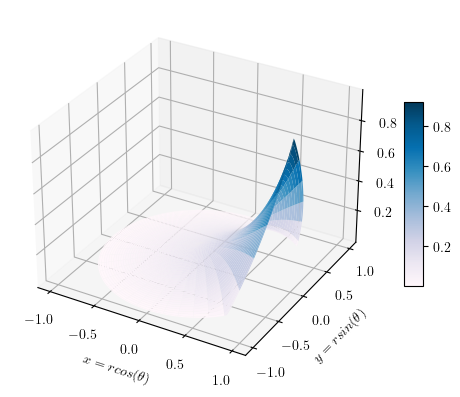}
        \hspace{0.01\linewidth}
        \includegraphics[width=0.48\textwidth]{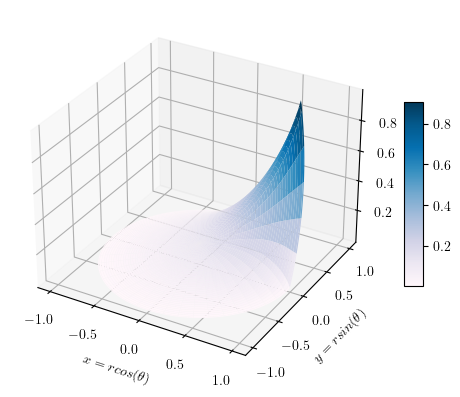}
        \caption{Basis functions $\phi_0$ and $\phi_1$ for $\omega = \pi/4$.\label{fig:sd_basis_2d}}
\end{figure}

Ultimately, both Green's functions method and the method of fundamental solutions were used to generate the basis functions. Both methods produced nearly identical results, with differences in the order of $1e{-6}$ near the boundary. Figure \ref{fig:sd_basis_2d} shows the shape of the basis functions within the domain~$K$ with radius $a=1$ and $\omega = \pi/4$.

We proceed to construct the space $\Xi_\varepsilon (K)$ recursively as detailed in the proof of \Cref{thm:one}. In the examples presented, the domains $K$ and $D$ are chosen as concentric disks centred at the origin with radii $0.5$ and $3$, respectively. The construction process begins with selecting the number of layers~$L$, followed by determining the number $n$ of basis functions per layer; note that $n$ is chosen to be constant across all layers. In \Cref{thm:one}, the number of basis functions on each layer is given by $n \sim c L^{d-1}$, where $c$ is a constant that depends on the distance from $K$ to the boundary $\partial D$, the constant $c_R$ in \eqref{eq:innreg}, and the approximation properties of the space $V_n(K)$. For our experiments, this constant is assumed to be $c = 2$ so that $n = 2L$. The dimension of the approximation space~$\Xi_\varepsilon(K)$ is then given in terms of the number of layers $$n_\Xi:=\dim \Xi_\varepsilon(K) = L n = 2 L^2.$$
The following is an overview of the algorithm that calculates the $L^2$-approximation $\xi_u\in\Xi_\eps(K)$ of~$u\in X_{0,0}(D)$. Since the support of the basis functions $\phi_i$ extends over the entire domain~$K$, the mass matrix $M \in \R^{n\times n}$ is not sparse. Nevertheless, the basis functions $\phi_i$ differ only by an angular shift, as shown in Figure \ref{fig:sd_basis_2d}, with $\phi_{i+1}(r,\theta) = \phi_{i}(r,\theta-\omega)$. This property makes the mass matrix a circulant matrix, meaning $m_{ij} = m_{\abs{i-(j\bmod n)}}$ for $i,j = 1,\ldots,n$. Additionally, $M$ is also symmetric. Consequently, only $\left \lfloor \frac{n}{2} \right \rfloor +1$ unique entries need to be computed on each layer. The entries $m_{ij}$ were computed using adaptive quadrature rules namely the Gauss-Kronrod rules~\cite{kronrod}.\\
It is also worth mentioning that each basis function on a given layer $s$ is shifted by $\omega/2$ relative to its counterpart in the preceding layer $s-1$. This shift is introduced to reduce the linear dependence between basis functions in successive layers. Experimental results indicate that implementing this shift improves the approximation error.
\begin{algorithm}
\caption{Calculating the $L^2$-approximation of a harmonic function}
    \begin{algorithmic}
        \State $s\gets 1$
        \While{$s\leq L$}
        \State generate the piecewise linear functions $\psi_i^{(s)}$, $i=1,\ldots,n$, on the boundary of the domain $K_s$
        \State calculate the basis representation $\tilde \phi^{(s)}_i$, $i=1,\ldots,n$, using MFS or Green's functions
        \State assemble the mass matrix $M^{(s)}$ with entries $m^{(s)}_{ij} = (\tilde \phi^{(s)}_{i}, \tilde \phi^{(s)}_{j})_{L^2}$ 
        \State calculate the right-hand side $b^{(s)}$, where $b^{(s)}_{i} = (u - \sum_{l=1}^{s-1} c^{(l)} \cdot \tilde \phi^{(l)}, \tilde \phi^{(s)}_i)_{L^2(K_s)}.$
        \State solve the system $M^{(s)} c^{(s)} = b^{(s)}$.
        \EndWhile
        \State compute $\xi_u := \sum_{s=1}^L c^{(s)} \cdot \tilde \phi^{(s)}|_{K_s}$
    \end{algorithmic}
\end{algorithm}

\subsection*{Trigonometric Approximation}
Instead of the piecewise linear approximation on the boundary of the disk, we can use trigonometric approximation.
Any $2\pi$-periodic $C^r$-function $f:[0,2\pi)\to\R$ can be approximated by a trigonometric polynomial $s_n[f]\in T_n:=\{t_n(\theta)=a_0+\sum_{j=1}^{n-1} a_j\cos(j\theta)+b_j\sin(j\theta)\}$ such that
\[
\norm{f-s_n[f]}_{L^\infty[0,2\pi)}\lesssim n^{-r};
\]
see \cite{MR36:571}. Define $V^T_n(K):=\{v\in X_{0,0}(K) \text{ such that } \gamma v\in T_n \}$. It can be easily seen that
the elements of $V_n^T(K)$ have the form
\[
v(r,\theta)=a_0+\sum_{j=1}^{n-1} \left(\frac{r}{a}\right)^j[a_j\cos(j\theta)+b_j\sin(j\theta)].
\]
Since $u\in X_{0,0}(D)$ is $C^\infty$ on $\partial K$, we obtain from the continuity of the trace operator
\[
\min_{v\in V^T_n(K)}\norm{u-v}_{H^{1}(K)}\lesssim \min_{v\in T_n}\norm{u-v}_{H^{1/2}(\partial K)}\lesssim \norm{u-s_n[u]}_{L^\infty(\partial K)}
\lesssim n^{-r}
\]
for any $r\in\N$. Hence, $V^T_n(K)$ can benefit from any order of regularity.

\subsection*{Examples}
We present the results and convergence error analysis for two harmonic functions. The approximation was computed for various numbers of layers $L$ and numbers of basis functions $n$, adhering to the relation $n = 2L$. For each scenario, the basis functions were computed using the method of fundamental solutions and Green's functions. In the MFS, the number of singularities was fixed at $N=256$, uniformly distributed on a circle of radius $R_s=a_s+0.01$, where $a_s$ represents the radius of $K_s$. According to \Cref{cor:erreps} the error estimate becomes $\norm{u-\xi_u}_{L^2(K)} \leq \eps \norm{u}_{L^2(D)}$, where $\varepsilon \lesssim \exp(-2 \sqrt{\dim \Xi_\varepsilon (K)})$.
% implying that the error behaves as $\varepsilon \lesssim \exp(-\sqrt{\dim \Xi_\varepsilon (K)})$.

\begin{figure}[h!t]
     \centering
     \subfigure[$\abs{u_1 -\xi_{u_1}}$ for $8$ layers with $16$ functions per layer.]{
        \includegraphics[width=0.45\linewidth]{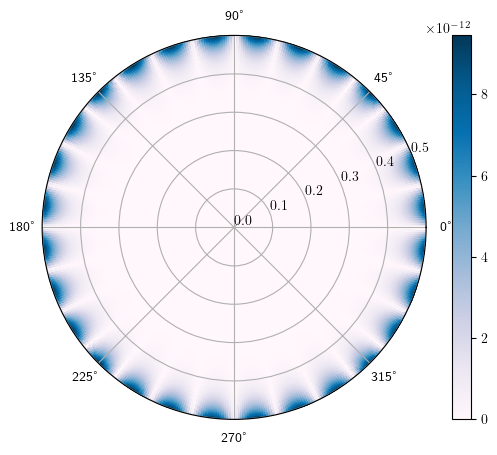}}
     \hspace{0.01\linewidth}
     \subfigure[$L^2$-error versus the number of degrees of freedom.]{     
         \includegraphics[width=0.5\linewidth]{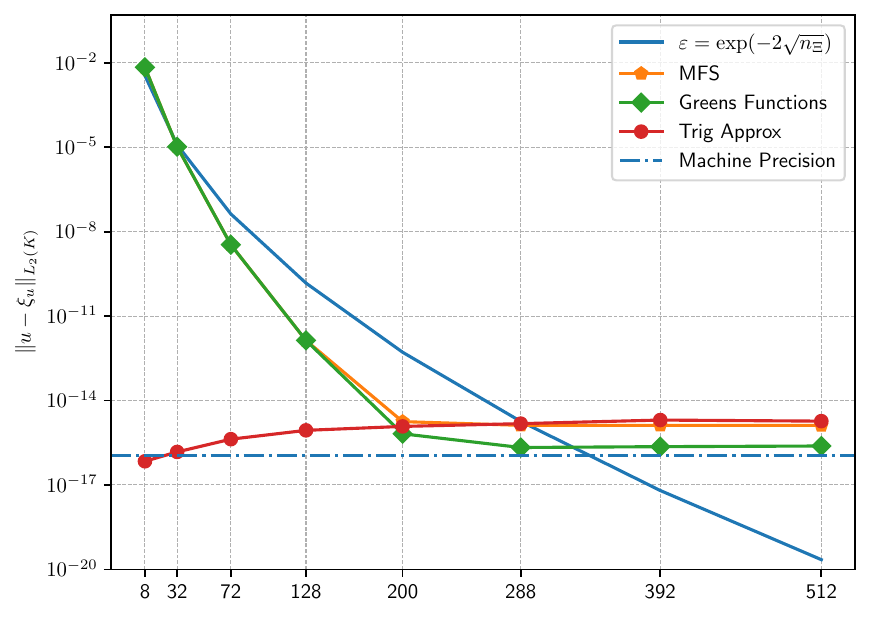}}
    \caption{Error and order of convergence for the harmonic function $u_1(r,\theta) = r^2\sin(2\theta)$.}
    \label{fig:sd_exp_conv_1}
\end{figure}

\begin{figure}[ht]
    \centering
     \subfigure[$\abs{u_1 -\xi_{u_1}}$ for $8$ layers with $16$ functions per layer.]{
        \includegraphics[width=0.45\linewidth]{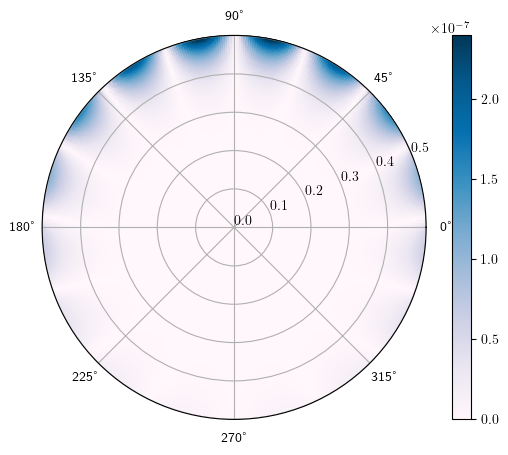}}
     \hspace{0.01\linewidth}
     \subfigure[$L^2$-error versus the number of degrees of freedom.]{     
         \includegraphics[width=0.5\linewidth]{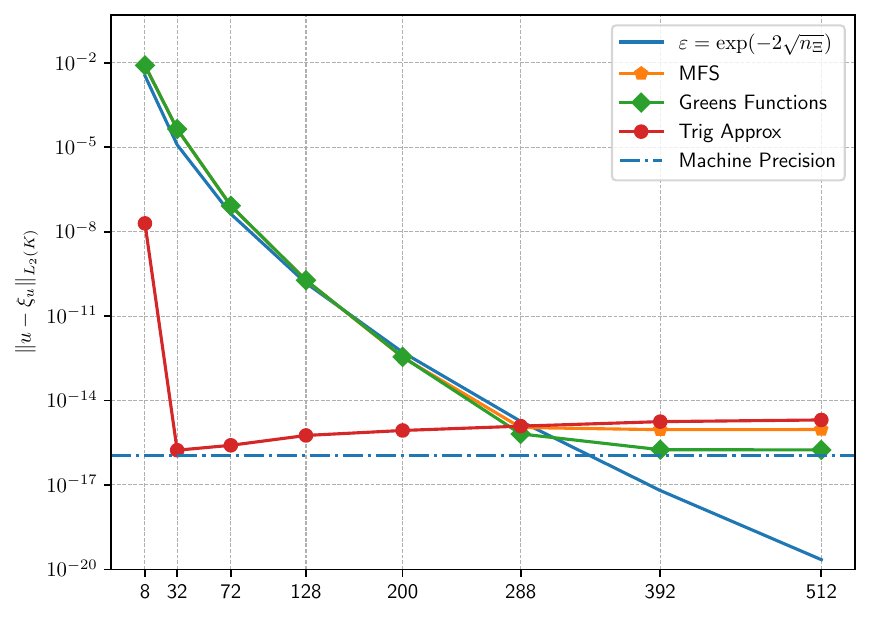}}
        \caption{Error and order of convergence for ${u_2(r,\theta)=\exp(r\sin(\theta)) \sin(r\cos(\theta))}$.}
        \label{fig:sd_exp_conv_2}
\end{figure}

Figures \ref{fig:sd_exp_conv_1} and \ref{fig:sd_exp_conv_2} show the error behaviour for two harmonic functions $u_1(r,\theta) = r^2 \sin(2\theta)$ and $u_2(r,\theta) = \exp(r\sin(\theta)) \sin(r\cos(\theta))$, respectively. Subfigure~(a) shows the pointwise absolute error~$\abs{u-\xi_u}$ using MFS with $8$ layers and $16$ basis functions per layer. For the second function, the absolute error is higher for the same number of basis functions due to the exponential term, which introduces unbounded derivatives near the upper part of the boundary; see Figure~\ref{fig:sd_exp_conv_2}.

Subfigure~(b) plots the approximation error in the $L^2$-norm versus the number of degrees of freedom, i.e.\ the dimension of $\Xi_\eps(K)$. Both functions exhibit exponential convergence, faster than or similar to the theoretical prediction, with the error saturating at numerical zero, approximately $1e{-16}$, due to machine precision. MFS and Green's functions yield similar results, though the latter achieves slightly lower approximation errors for a larger number of basis functions. It can also be seen that the convergence rate for the second function is slower than that of the first, which is again due to the exponential term in the function.

An additional test was conducted for both functions, utilizing the trigonometric approximation method described in the previous section. In this context, the number of degrees of freedom corresponds to the number of frequencies $n$ used to approximate the function. The results indicate that for the first harmonic function, $u(r,\theta) = r^2\sin(2\theta)$, the error remains at numerical zero for all $n\geq 2$, as the function can be exactly represented using only two frequencies. For the second harmonic function, the approximation exhibits a convergence rate that surpasses any algebraic order, which is consistent with theoretical predictions for infinitely smooth functions. However, this method is limited to highly smooth domains and cannot be readily extended to the general case of Lipschitz domains.

\subsection{Lipschitz Domains}
In this section, we present the numerical construction of exponentially convergent spaces~$\Xi_\eps(K)$ on Lipschitz domains. The steps are similar to those presented in the previous sections. Once again we start by constructing the basis of the space $\op S^1_h(\partial K)$, however now $K \subset \R^2$ is a square centred at the origin with a side length of $2a$. First, the boundary of the square is parametrised in terms of a free parameter $t$ as follows
\[
    [0,4) \ni t \mapsto \begin{bmatrix}x(t)\\y(t)\end{bmatrix} = 2a\begin{bmatrix}\gamma \left(t-\frac{3}{2}\right)\\\gamma \left(t-\frac{5}{2} \right)\end{bmatrix},\]
where $\gamma(z) = \max \left(-\frac{1}{2},\min\left(\frac{1}{2},1-\abs{z}\right) \right)$. 
Then, the boundary mesh $\op T_h$ on $\partial K$ is defined as
\begin{equation*}
    \op T_h := \left \{ \left [ t_j,t_{j+1} \right] : 0\leq j \leq n-1 \right \} ,
\end{equation*}
where $0=t_0 < \ldots < t_n = 4$ with uniform mesh size $h:=8a/n$. It is obvious that the mesh $\op T_h$ is periodic, meaning $t_0 = t_n$, so the mesh has exactly $n$ intervals and $n$ distinct nodes. 

\begin{figure}[ht]
    \centering
    \includegraphics[width=0.48\linewidth]{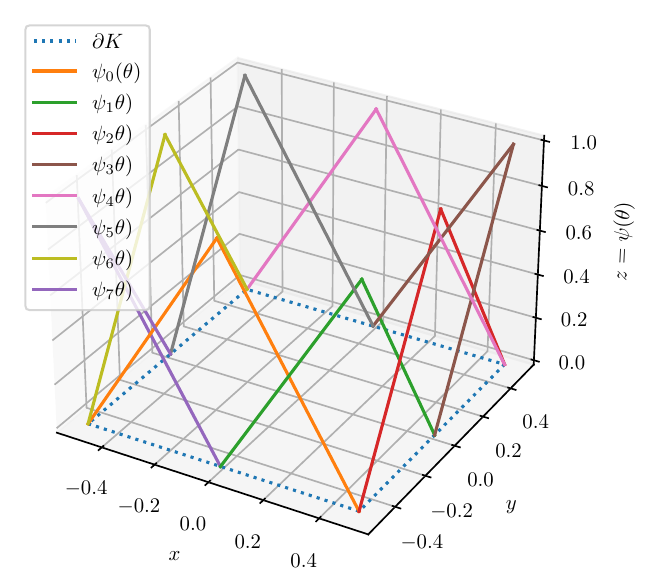}
    \caption{Piecewise linear functions $\psi_i$ on $\partial K$ for $\omega = 1/2$. \label{fig:ld_hf}}
\end{figure}

Let $\omega := h /(2 a)$ be the discretisation parameter of $t$, we define the basis $\psi_i$ as
\[
\psi_i(x(t),y(t)) = 
  \begin{cases}
    \frac{t-(i-1)\omega}{\omega}, &  (i-1)\omega \leq t \le i\omega, \\
    \frac{(i+1)\omega -t}{\omega}, & i\omega \leq t \le (i+1)\omega,\\
    0, &\textnormal{otherwise.}
  \end{cases}
\]
The shapes of the basis functions $\psi_i$ for different discretisation parameters $\omega$ are shown in Figure~\ref{fig:ld_hf}.

\begin{figure}[ht]
     \centering
     \includegraphics[width=0.48\linewidth]{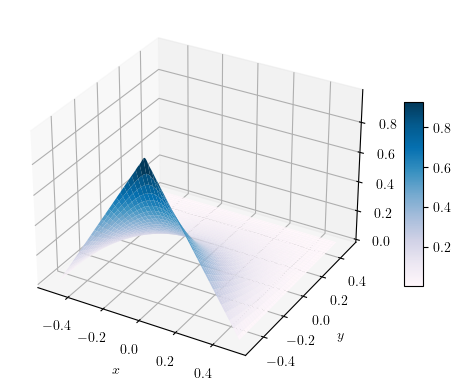}
        \hspace{0.01\linewidth}
     \includegraphics[width=0.48\textwidth]{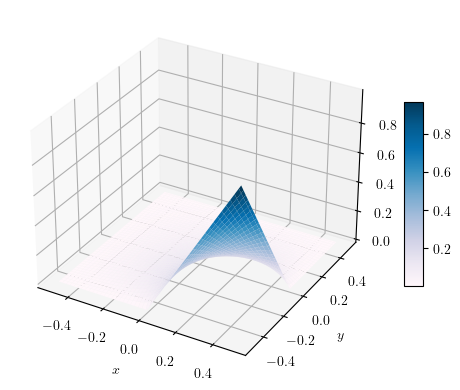}
        \caption{Basis functions $\phi_1$ and $\phi_2$ for $\omega = 1/2$.\label{fig:ld_basis}}
\end{figure}

In a similar manner to the previous section, we compute the basis functions $\phi_i$, $i=1,\dots,n$, of the space $V_n(K)$ as the harmonic extension of the basis functions $\psi_i$ by solving the boundary value problems
\begin{subequations}    \label{eq:ld_harm_bas}
    \begin{alignat}{3} 
        -\Delta \phi_i(x,y) &= 0 &\quad& \text{in } K ,\\
        \phi_i(x(t),y(t)) &= \psi_i(t) &\quad& \text{on } \partial K.
    \end{alignat}	
\end{subequations}
In the case of smooth domains, both Green's functions method and the method of fundamental solutions were considered, however, Green's functions are not known for general Lipschitz domains with inhomogeneous boundary conditions. Instead, the finite element method provides a straightforward approach for generating the basis functions. Consequently, both the method of fundamental solutions and the finite element method were employed to solve \eqref{eq:ld_harm_bas}.\\
The main procedure of the method of fundamental solutions remains the same. The solution is represented as a linear combination of $N$ fundamental solutions with singularities positioned outside of the domain while the boundary of $K$ is discretised into $M$ collocation points to enforce the boundary conditions. For a square domain with side length $2a$, the singularities are placed on a square with side length $2(a+\delta a)$. Again, the number of collocation points $M$ is selected to be at a higher resolution than that of the singularities, namely $M\approx 4N$. This arrangement ensured proper alignment and accurate approximation, particularly at the corners of the domain. As for the finite elements method, a uniform grid was used with Courant elements. The stiffness matrix is assembled only once, as subsequent computations involve modifications solely in the right-hand side vector.
\begin{figure}[ht]
\centering
\begin{tikzpicture}[scale=1.8]
    \draw (1,1) -- (1,-1) -- (-1,-1) -- (-1,1) -- cycle;
    \node at (-0.5,-0.5) {$K$};
    \def\n{12}
    \def\step{2.2/\n}
    \foreach \i in {0,...,\n} {
        \filldraw[blue] ({-1.1 + \i*\step}, 1.1) circle (0.02);
        \filldraw[blue] ({-1.1 + \i*\step}, -1.1) circle (0.02);
        \filldraw[blue] (-1.1, {-1.1 + \i*\step}) circle (0.02);
        \filldraw[blue] (1.1, {-1.1 + \i*\step}) circle (0.02);
    }
    \draw[dashed] (0,0) -- (0,1);
    \node at (-1.1,1.1) [left,blue] {$q_j$};
    \draw (0,0) -- (1.1,0);
    \filldraw  (0,0) circle (0.03);
    \node at (0.1,0.5) {$a$};
    \node at (0.5,0.1) {$a+\delta a$};
\end{tikzpicture}
\captionof{figure}{Singularities placement outside of a square domain.}
\end{figure}
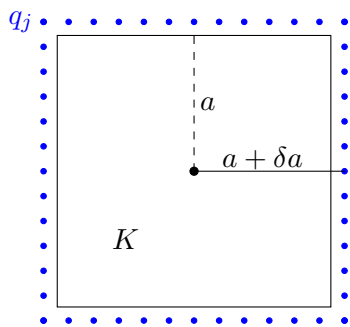
Figure \ref{fig:ld_basis} displays the shape of the basis functions within the domain $K$ with side length $2a=1$ and $\omega = 1/2$.

The space $\Xi_\varepsilon (K)$ is constructed recursively using the same procedure from the previous section. In the examples provided, the domains $K$ and $D$ are chosen as squares centred at the origin with side lengths $0.5$ and $3$, respectively. Also, the number of the basis on each layer $n$ was chosen again to be double the total number of layers $L$. In the previous section, we were able to exploit the rotational invariance of the basis functions $\phi_i$ to lower the numbers of elements $m_{ij}$ to be computed. In the case of Lipschitz domains, due to the presence of corners in the square domain, we cannot use such rotational invariance. In this case, only the symmetry of the $L^2$ inner product is exploited, making the number of unique matrix elements that needed to be computed approximately $n^2/2$ for each layer. It is also worth noting that, for Lipschitz domains, Clenshaw–Curtis quadrature rules \cite{clenshaw_curtis} were used instead of Gauss–Kronrod quadrature rules. The former proved to handle the integration more effectively in this case. Similar to the case of smooth domains, the basis functions on each layer were shifted by $\omega/2$ relative to the preceding layer.

\subsection*{Examples}
 We present the results and convergence error analysis for two harmonic functions. In the method of fundamental solutions, the number of singularities was fixed at $N=256$, uniformly distributed on a virtual square of side length $2(a_s +0.01)$, where $2 a_s$ is the side length of $K_s$. As for the finite element method, two discretisation parameters were tested: $N_\textnormal{FEM} = 512$ and $N_\textnormal{FEM} = 1024$, where $N_\textnormal{FEM}$ denotes the number of elements in one dimension.
\begin{figure}[ht]
     \centering
     \subfigure[$\abs{u_1 -\xi_{u_1}}$ for $8$ layers with $16$ functions per layer.]{
        \includegraphics[width=0.45\linewidth]{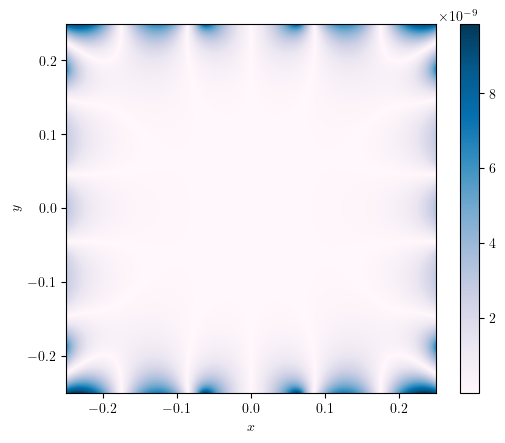}}
     \hspace{0.01\linewidth}
     \subfigure[$L^2$-error versus the number of degrees of freedom.]{     
         \includegraphics[width=0.5\linewidth]{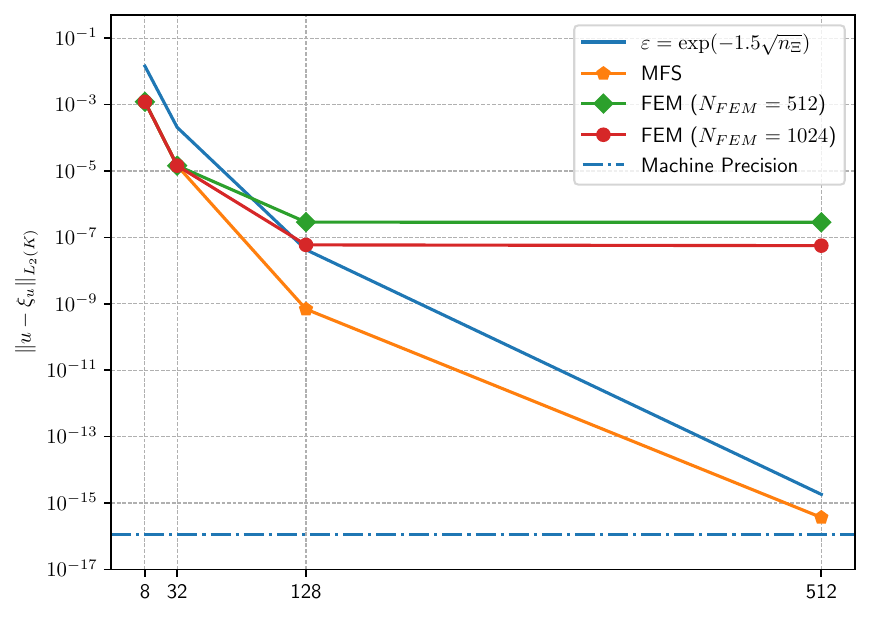}}
    \caption{Error and order of convergence for the harmonic function $u_1(x,y) = x^3-3xy^2$.\label{fig:ld_exp_conv_1}}

\end{figure}
\begin{figure}[h!t]
    \centering
     \subfigure[$\abs{u_1 -\xi_{u_1}}$ for $8$ layers with $16$ functions per layer.]{
        \includegraphics[width=0.45\linewidth]{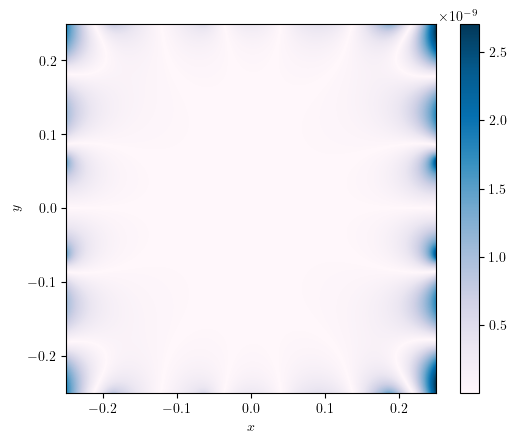}}
     \hspace{0.01\linewidth}
     \subfigure[$L^2$-error versus the number of degrees of freedom.]{ 
         \includegraphics[width=0.5\linewidth]{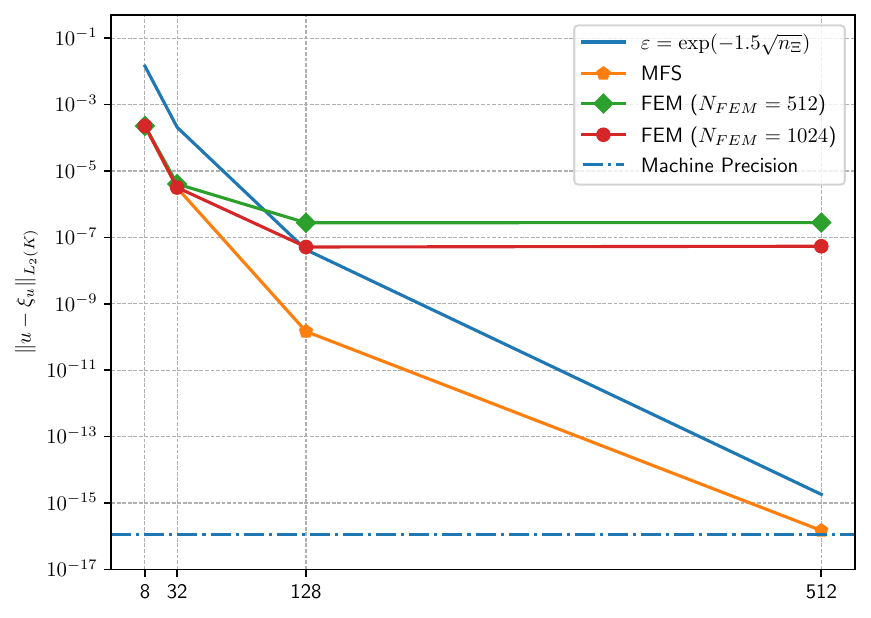}}
        \caption{Error and order of convergence for the harmonic function $u_2(x,y) = \exp(x) \sin(y) $. \label{fig:ld_exp_conv_2}}
\end{figure}
From Figures~\ref{fig:ld_exp_conv_1} and \ref{fig:ld_exp_conv_2} 
it can be observed that the numerical results for the method of fundamental solutions exhibit a slightly faster convergence rate than the one predicted by the theoretical analysis. In contrast, for the finite element method, the error shows exponential convergence for a small number of basis functions but saturates around $10^{-7}$ for both functions. This saturation can be attributed to the discrete harmonicity of the basis functions generated using FEM. Additionally, it is observed that the error saturates at a lower value as the number of elements $N_\textnormal{FEM}$ increases. Although finer approximations could be tested, this approach proved to be both time-consuming and memory-intensive, rendering it impractical.
\bibliographystyle{plain}
\bibliography{citations,Bib/analysis,Bib/hmatrix,Bib/buecher,Bib/SparseLU,Bib/fem_bem,Bib/multilevel,Bib/ApproxInterp}

\end{document}